\documentclass[10pt,a4paper,onecolumn]{article}
\usepackage{blindtext,amssymb}
\usepackage[utf8]{inputenc}
\usepackage[english]{babel}
\usepackage{graphicx}
\graphicspath{ {/} }
\usepackage{amsmath}
\usepackage{amsfonts}
\usepackage{ragged2e}
\usepackage{wrapfig}
\usepackage{a4wide}
\usepackage{epstopdf}
\usepackage{float}
\usepackage{caption}

\DeclareCaptionFont{mysize}{\fontsize{9}{9.6}\selectfont}
\usepackage{hyperref}
\newif\ifhyper
\hypertrue
\ifhyper
\hypersetup{
   citecolor = {green},
   colorlinks = {true}, % false
   urlcolor = {blue} % magenta
}
\fi

\newcommand{\rom}[1]{{{\color{black}{#1}}}}
\newcommand{\hlt}[1]{{{\color{black}{#1}}}}

\usepackage{listings}
\usepackage{color}
\usepackage{soul}

\definecolor{dkgreen}{rgb}{0,0.6,0}
\definecolor{gray}{rgb}{0.5,0.5,0.5}
\definecolor{mauve}{rgb}{0.58,0,0.82}

\begin{document}

\title{Analyzing Riemann's hypothesis}
\author{Mercedes Or\'us--Lacort$^{1,2}$, Rom\'an Or\'us$^{3, 4, 5}$, Christophe Jouis$^{6}$  \\ \\ \\
\multicolumn{1}{p{0.9\textwidth}}{\centering\emph{{
$^{1}$College Mathematics, Universitat Oberta de Catalunya, \\ Rambla del Poblenou 156, 08018 Barcelona, Spain \\ 
$^{2}$ College Mathematics, Universitat Nacional de Educacion a Distancia, \\ Calle Pintor Sorolla 21, 46002, Valencia, Spain \\
$^{3}$ Donostia International Physics Center, \\ Paseo Manuel de Lardizabal 4, E-20018 San Sebasti\'an, Spain \\
$^{4}$ Multiverse Computing, \\ Paseo de Miram\'on 170, E-20014 San Sebasti\'an, Spain \\
$^{5}$ Ikerbasque Foundation for Science, \\ Maria Diaz de Haro 3, E-48013 Bilbao, Spain \\
$^{6}$ Laboratoire d'informatique de l'\'Ecole polytechnique, LIX, UMR 7161, 1 rue Honor\'e d'Estienne d'Orves, B\^atiment Alan Turing, Campus de l'\'Ecole Polytechnique, 91120 Palaiseau, France}}}}

\maketitle

\begin{abstract}
In this paper we perform a detailed analysis of Riemann's hypothesis, dealing with the zeros of the analytically-extended zeta function. We use the functional equation $\zeta(s) = 2^{s}\pi^{s-1}\sin{(\displaystyle \pi s/2)}\Gamma(1-s)\zeta(1-s)$ for complex numbers $s$ such that $0<{\rm Re(s)}<1$ and the reduction to the absurd method where we use an analytical study based on a complex function and its modulus as a real function of two real variables in combination with a deep numerical analysis to show that the real part of the non-trivial zeros of the Riemann zeta function is equal to $1/2$ to the best of our resources. This is done in two steps. Firstly,  we show what would happen if we assumed that the real part of $s$ has a value between $0$ and $1$ but different from $1/2$ arriving at a possible contradiction for the zeros. Secondly assuming that there is no real value $y$ such that $\zeta\left(1/2 +yi \right)=0$ by applying the rules of logic to negate a quantifier and the corresponding Morgan's law we also arrive to a plausible contradiction. Finally, we analyze what conditions should be satisfied by $y \in \mathbb R$ such that $\zeta(\displaystyle 1/2 +yi)=0$. While these results are valid to the best of our numerical calculations, we do not observe and foresee any tendency for a change. Our findings open the way towards assessing the validity of Riemman's hypothesis from a fresh and new mathematical perspective.  

\end{abstract}


\begin{ams}
11--02
\end{ams}

\section{Introduction}
\label{sec1} 

Riemann's hypothesis, first formulated by Bernhard Riemann in 1859 \cite{4}, is a conjecture about the distribution of the zeros of the Riemann zeta function $\zeta(s)$ \cite{1}. Due to its relationship with the distribution of prime numbers in the set of natural numbers, proving this hypothesis is one of the most important open problems in contemporary mathematics \cite{2,3}. 

\hlt{In this paper, we analyze Riemman's hypothesis by dealing with the zeros of the analytically-extended zeta function. To be specfic, we make use of the functional equation $\zeta(s) = 2^{s}\pi^{s-1}\sin{(\displaystyle \pi s/2)}\Gamma(1-s)\zeta(1-s)$ for complex numbers $s$ such that $0<{\rm Re}(s)<1$. Our objective is to assess the possibility that the only non-trivial zeros of the Riemann zeta function are those complex numbers whose real part is equal to $1/2$. We pursue this analysis using a combination of analytical and numerical techniques, as discussed in detail below.} 

\hlt{
\subsection{Previous research}}

Riemann hypothesis has recently been the subject of much research. For instance,  \cite{r1} establishes an analogue to Lagaria's criterion for the hypothesis in terms of harmonic series.  \cite{r2} describes the connection between the zeta function and the solution of the Majorana fermionic equation in curved space-times.  \cite{r3} studies the validity of the hypothesis via an statistical analysis of Mertens function. In addition,  \cite{r4} makes use of Riemann's hypothesis to obtain asymptotic formulas for the second moment of the n$th$ antiderivative of the argument of the zeta function.  \cite{r5} studies ergodic theorems to provide new characterisations of the hypothesis. In  \cite{r6} the authors observe that the zeta function corresponds to the spectrum of a certain quantum Hamiltonian capturing the near-horizon dynamics of the Schwarzschild black hole.  \cite{r7} describes an algorithm to compute very high Riemann zeros using random walks. The work in  \cite{r8} provides a Riesz-type criterio for the generalized Riemann hypothesis. In  \cite{r9}, the authors consider and analogue of the hypothesis and quantum walks. A Riesz-type criteria for the hypothesis in terms of one variable was also provided in  \cite{r10}.  \cite{r11} considers deformations of the Keiper-Li sequence to analyze the hypothesis. The authors of  \cite{r12} propose a positivity conjecture for matrices related directly to the hypothesis.  \cite{r13} studies the meromorphic extensions of fractal zeta functions from quasiperiodic sets. Also regarding meromorphic functions,   \cite{r14} considers the meromorphic modular forms of a family of generalized $L$-functions, and relates them to the zeta function. 

Additionally, other important research works were published recently that have a connection to our work. An example are the new results by Conrey \cite{r15}, who analyzes that the major difficulty when trying to prove the hypothesis through analysis come from the fact that the zeros of $L$-functions have a very different behavior to the zeros of many of the special mathematical functions  and only recently it was found that the modularity of the $L$-functions is associated to elliptic curves, which could help in this direction. Moreover, recent work by Liu \cite{r16} attempts to prove the Riemann hypothesis for both the Riemann zeta-function $\zeta(s)$ and the Dirichlet $L$-function $L(s,\chi)$   through an equivalent condition on the Farey series set forth by Franel and Landau. Finally, Liu and Wang analyzed recently \cite{r17} the Riemann problem of the high-order Jaulent-Miodek (JM) equation with initial data of step discontinuity, as explored by Whitham modulation theory, and found that the periodic wave solutions of the high-order JM equation are described by the elliptic function along with the Whitham modulation equations. 

\hlt{The hypothesis has also been the subject of intensive research works in the past in several ambits, see for instance \cite{n1, n2, n3, n4, n5, n6, n7, n8, n9, n10} and references therein.}

\hlt{
\subsection{Everyday life applications}
Even if belonging to the field of pure mathematics, this hypothesis also has applications in our daily life. We have already mentioned some of them above, when enumerating references to previous important and recent research works. To put everything in context, let us discuss briefly some of these applications here.

The hypothesis also finds a wide spectrum of applications in science and technology. For instance, the zeros of the zeta function have important connections to the energy spectrum of classical chaotic systems \cite{chaos}, quantum Hamiltonians \cite{german}, as well as to scattering amplitudes in quantum field theory \cite{qft}. Quantum physics is currently about us, and in the near future computers will process information directly at a quantum level. The hypothesis has therefore a direct impact on the upcoming quantum technologies and their industrial implications \cite{quantum}.

In addition, there are also important implications in the field of cryptography. The zeros of the zeta function can be interpreted as harmonic frequencies in the distribution of primes, leading to studies of the distribution of distances between consecutive primes \cite{cramers}. Such analysis is key in the security of asymmetric-key cryptosystems such as RSA, which is based on the fact that finding the prime factorization of a natural number is a hard computational problem \cite{RSA}, though not for quantum computers \cite{shor}. Though RSA is no longer used as an standard cryptographic protocol (in favor of symmetric-key schemes), it has been a technical standard for many years and still, as of today, is used in certain non-critical applications.

There are many other daily life applications of the hypothesis. Apart from encryption algorithms, which rely heavily on the properties of prime number distributions, zeros of the zeta function are also related to particle distributions in quantum statistical mechanics, as well as to the eigenvalues of random matrices \cite{randommatrix}. As such, random matrices have by themselves a very wide spectrum of practical applications, including the energies of heavy uranium-like nucleai, the behaviour and dynamics of financial markets, and even the development of new machine learning and Artificial Intelligence (AI) algorithms based on neural networks and deep learning. This last point is intriguing, since new techniques are required to boost performance and decrease the energy consumption of current AI models, such as Large Language Models \cite{LLM}. In this overall context, the hypothesis impacts also on current advanced software used for cybersecurity and artificial intelligence, such as the EC3 Software \cite{EC3}.} 

\hlt{
\subsection{About this paper}} 

\hlt{As stated above, here we analyze Riemman's hypothesis by dealing with the zeros of the analytically-extended zeta function using the functional equation $\zeta(s) = 2^{s}\pi^{s-1}\sin{(\displaystyle \pi s/2)}\Gamma(1-s)\zeta(1-s)$ for complex numbers $s$ such that $0<{\rm Re}(s)<1$.

The objective of this analysis is to assess the possibility that the only non-trivial zeros of the Riemann zeta function, using its functional equation, are those complex numbers whose real part is equal to $1/2$. More specifically, we want to show that if $x$ is a real number such that $0 < x < 1$, then it can be satisfied that there exists $y\in {\mathbb R}$ such that $\zeta(x+iy)=0$ if and only if $x=1/2$.

In order to develop our analysis, we will use a combination of two mathematical methods. The first method is an analytical technique called \emph{reduction to the absurd}, and the second one is based on \emph{numerical analysis} techniques using two different types of specialized mathematical software: Matlab and Wolfram Alpha.

The reduction to the absurd method is a very usual and well-known technique, when we want to prove that a conditional statement such as ``If A happens, then B happens" is true.

As we shall see, in this paper we make use that, if the statement ``If A happens, then B happens" is true, then by using the fundamental laws of mathematical logic, the negated statement must also be true. That is, ``If B does not happen, then A does not happen" must also be true.

Similarly here we may assume that ``A happens" but ``B does not happen". In this way if we are able to show that if ``A does not happen", then we arrive at what we call \emph{a contradiction}, because we were assuming that ``A happens". Therefore, the statement ``If B does not happen, then A does not happen" is true, and ``If A happens, then B happens" is also true.

Additionally, the numerical analysis method is based on algorithms developed using Matlab and  Wolfram Alpha. Note that, like all techniques based on numerical analysis, in our algorithms some of the calculations will be iterated until a sufficiently small error is attained, which we take as convergence criteria. In our numerical approach, both Matlab and Wolfram Alpha work up to a preset error given by machine precision. Therefore, of course, the numerical analysis methods work under certain assumptions, for example not exceeding a certain accepted (minimum and small enough) error, as well as other limitations coming from potential sources of error such as floating point errors. 

As for specific algorithms, in Matlab we will use the ``gradient-free Nelder-Mead algorithm" to minimize two functions to the best of our numerical capabilities. Wolfram Alpha is used to solve a non-linear system of equations. See Appendix \ref{apA} for more detailed information on our numerical methods.  

Diving into more details of our actual derivations, as a first step, we show two Propositions that we use later to work analytically throughout the whole paper.}

Using the reduction to the absurd method, we develop an analytical study based on a complex function, and its modulus as a real function of two real variables. And at a certain point in our study, we combine it with intensive numerical analysis at some steps. Our derivations are compatible with the real part of the non-trivial zeros of the Riemann zeta function being equal to $1/2$, to the best of our resources. 

We do this in two steps. First, we show what would happen if we assume that the real part of $s$ has a value between $0$ and $1$ but different from $1/2$, arriving to a possible contradiction for the zeros. Second, assuming that there is no real value $y$ such that $\zeta\left(1/2 +yi \right)=0$, and by applying the rules of logic to negate a quantifier together with the corresponding Morgan's law, we also arrive to a plausible contradiction.

Finally, we also analyze what conditions should be satisfied by $y \in \mathbb R$ such that $\zeta(\displaystyle 1/2 +yi)=0$. While most of our results are fully analytic, at some specific parts of the analysis we need to rely on heavy numerical calculations. Some part of our analysis is therefore dependent on them. However, we do not observe nor foresee any tendency for a change in our calculations, which leads us to conjecture that the validity of our conclusions is general. Our approach also opens the way towards a new mathematical concept to assess Riemman's hypothesis. 

The paper is organized as follows. In Sec.\ref{sec2} we review the basics of the hypothesis. In Sec.\ref{sec3} we implement our analysis, as briefly outlined above. In Sec.\ref{seczeros} we provide analytical constraints to be satisfied by the zeros that come out naturally from our analysis. Finally, in Sec.\ref{sec4} we summarize  our conclusions and perspectives for future work. In addition, in Appendix \ref{apA} we provide the Matlab computer codes used, and give a brief explanation of the numerical techniques used both with Matlab and with Wolfram Alpha. 

\section{Zeta function and Riemman's hypothesis} 
\label{sec2} 

The Riemann zeta function $\zeta(s)$ is defined in complex numbers as the sum of an infinite series as follows:
\begin{equation} \label{eq1}
\zeta(s)=\sum_{n=1}^{\infty}\displaystyle \frac{1}{n^{s}}. 
\end{equation}
The series is convergent when ${\rm Re}(s)$ is strictly greater than $1$. Leonhard Euler showed that this series is equivalent to Euler's product, 
\begin{equation} \label{eq2}
\zeta(s)=\prod_{p {\kern 1pt} {\kern 1pt} {\rm{ }}{\kern 1pt} {\rm prime}}\displaystyle \frac{1}{1-p^{-s}},
\end{equation}
where the infinite product extends over the set of all prime numbers $p$, and again converges for complex $s$ whose real part is greater than $1$. The convergence of the Euler product shows that $\zeta(s)$ has no zeros in this region, since none of the factors in the product have zeros. 

Riemann's hypothesis deals with the zeros outside the convergence radius of the series in Eq.(\ref{eq1}) and/or the Euler product in Eq.(\ref{eq2}). To preserve the meaning of this hypothesis, one needs to analytically continue the zeta function $\zeta(s)$, so that it makes sense for any value of $s$. Any choice of extension will lead to the same conclusions as above, since the zeta function is meromorphic.  Hence, in particular, for complex numbers $s$ such that $0<{\rm Re}(s)<1$, the function $\zeta(s)$ can be expressed by the functional equation
\begin{equation} \label{eq3}
\zeta(s)=2^{s}\pi^{s-1}\sin{\left(\displaystyle \frac{\pi s}{2}\right)}\Gamma(1-s)\zeta(1-s), 
\end{equation}
where $\Gamma(s)$ is the Gamma function, defined by 
\begin{equation}
\Gamma(s) = \int_0^\infty t^{s - 1} e^{-t} dt. 
\end{equation} 
In this work we will make extensive use of this functional form. 

Initially, we consider the zeros of the above analytical extension of the zeta function. Some of these zeros are called "trivial", since they can be easily seen by inspection. In particular, from Eq.(\ref{eq3}) one can see that $s$ = $-2, -4, -6, \cdots$, i.e., all negative even integers, are trivial zeros, since they cancel the trigonometric function. Likewise, there are other (complex) values of $s$ such that $0 <{\rm Re} (s) <1$ and for which the zeta function also vanishes, which are called "non-trivial" zeros. The Riemann conjecture refers specifically to these non trivial zeros, stating the following:
\begin{center}
\textit{The real part of all non-trivial zeros of the Riemann zeta function is equal to $\displaystyle \frac{1}{2}$.}
\end{center}
The conjecture therefore implies that all the non-trivial zeros should lie on the critical line $s =\displaystyle 1/2+it$, where $t$ is a real number and $i$ is the imaginary unit. 

\section{Analysis of the hypothesis} 
\label{sec3} 

Let us start our analysis with the following two propositions: 

\noindent
\underline{\textbf{Proposition 1}}: {\it Using Riemann's zeta functional equation, it is satisfied that, $\forall s \in \mathbb C$ such that $0<{\rm Re}(s)<1$ and $\zeta(s)=0$, then $\zeta(1-s)=0$.}

\noindent
\underline{Proof:} Given Riemann's zeta functional equation in Eq.(\ref{eq3}), then $\zeta(1-s)$ is equal to 
\begin{equation} \label{eq5}
\zeta(1-s)=2^{1-s}\pi^{-s}\sin{\left(\displaystyle \frac{\pi (1-s)}{2}\right)}\Gamma(s). \zeta(s).
\end{equation}
Therefore, if there exists an $s \in \mathbb C$ such that $0<{\rm Re}(s)<1$ and $\zeta(s)=0$, then by the previous equation we have that 
\begin{equation} \label{eq6}
\zeta(1-s)=2^{1-s}\pi^{-s}\sin{\left(\displaystyle \frac{\pi (1-s)}{2}\right)}\Gamma(s)\cdot 0 = 0, 
\end{equation}
which proves the proposition. $\square$ 

\bigskip 

\noindent
\underline{\textbf{Proposition 2}}: {\it It is satisfied that $\zeta(a-bi)=\overline{\zeta(a+bi)}$.}

\noindent
\underline{Proof:} Using the definition of $\zeta(a+bi)$, we have that 
\begin{equation} \label{eq7}
\zeta(a+bi)=\sum_{n=1}^{\infty}\displaystyle \frac{1}{n^{a+bi}}. 
\end{equation}
This implies that 
\begin{eqnarray} 
\zeta(a+bi)&=&\sum_{n=1}^{\infty}\displaystyle \frac{n^{-bi}}{n^{a}} =\sum_{n=1}^{\infty}\displaystyle \frac{e^{\ln{(n^{-bi})}}}{n^{a}},  \nonumber \\ 
&=&\sum_{n=1}^{\infty}\displaystyle \frac{e^{-ib \ln{(n)}}}{n^{a}} = \sum_{n=1}^{\infty}\displaystyle \frac{\cos{(b \ln{(n)})}-i\sin{(b \ln{(n)})}}{n^{a}}, \nonumber \\
&=& \sum_{n=1}^{\infty}\displaystyle \frac{\cos{(b \ln{(n)})}}{n^{a}} - i\sum_{n=1}^{\infty}\displaystyle \frac{\sin{(b \ln{(n)})}}{n^{a}}. 
\end{eqnarray}
If we do the same but for $\zeta(a-bi)$, we have that 
\begin{eqnarray} 
\zeta(a-bi)&=&\sum_{n=1}^{\infty}\displaystyle \frac{n^{bi}}{n^{a}} =\sum_{n=1}^{\infty}\displaystyle \frac{e^{\ln{(n^{bi})}}}{n^{a}},  \nonumber \\ 
&=&\sum_{n=1}^{\infty}\displaystyle \frac{e^{ib \ln{(n)}}}{n^{a}} = \sum_{n=1}^{\infty}\displaystyle \frac{\cos{(b \ln{(n)})}+i\sin{(b \ln{(n)})}}{n^{a}}, \nonumber \\
&=& \sum_{n=1}^{\infty}\displaystyle \frac{\cos{(b \ln{(n)})}}{n^{a}} + i\sum_{n=1}^{\infty}\displaystyle \frac{\sin{(b \ln{(n)})}}{n^{a}}.
\end{eqnarray}
Therefore, comparing the above equations, we see that $\zeta(a-bi)=\overline{\zeta(a+bi)}$, as we want to prove. $\square$ 

We now go into the actual details of Riemann's conjecture. As a reminder, this reads as follows: 

\noindent
\underline{\textbf{Riemann hypothesis: }} {\it Let $x$ be a real number such that $0<x<1$. Then, it is satisfied that $\exists y \in \mathbb R$ such that $\zeta(x+yi) = 0$ if and only if $x=\displaystyle 1/2$.}

Our approach to the conjecture is as follows. Let $x$ be a real number such that $0<x<1$. We assess the validity of the double implication, and start with the first direction. 

\noindent
\underline{\textit{(i) If $\exists y \in \mathbb R$ such that $\zeta(x+yi) = 0$ then $x=\displaystyle 1/2$. }}

Our strategy is to use the reduction to the absurd method to validate it. With this in mind, let us then assume that there exists $s=a+bi \in \mathbb C$ such that $a\neq  \displaystyle 1/2$, $0 < a <1$ and such that $\zeta (a+bi) = 0$. 

\rom{If the above were true, then using Proposition 1 it would also be satisfied that $\zeta (1-(a+bi)) = 0$. Therefore, since the function $\zeta (s)$ is continuous $\forall s \in \mathbb C$ such that ${\rm Re}(s)\in(0, 1)$, then if $a$ is, for example, less than $1- a$, then this implies that $|\zeta (a+bi)| = |\zeta (1-(a+bi))| = 0$, and then also that $|\zeta (a+bi)|^2 = |\zeta (1-(a+bi))|^2 = 0$. 

Then, $|\zeta (x+yi)|^2$, as a real function of real variables $x$ and $y$, should have a maximum or a minimum for some $x + yi \in \mathbb{C}$, where $0 < a < x < 1-a < 1$ (assuming, without loss of generality, that $a < 1-a$). 

Let us now analyze if this is possible. First, notice that 
\begin{equation}
|\zeta (x+yi)|^2 = \zeta (x+yi) \overline{\zeta (x+yi)} = \zeta (x+yi) \zeta (x-yi), 
\end{equation}
where previously we used Proposition 2. Next, we calculate the following first-order partial derivative with respect to $x$: 
\begin{eqnarray}
\partial_x |\zeta (x+yi)|^2 &=& \partial_x \left(\zeta (x+yi) \zeta (x-yi) \right) = \left(\partial_x\zeta (x+yi)\right) \zeta (x-yi) + \zeta (x+yi) \left( \partial_x\zeta (x-yi) \right), \nonumber \\  
&=& \zeta' (x+yi) \zeta (x-yi) + \zeta (x+yi) \zeta' (x-yi) = 0. 
\end{eqnarray}
This equation implies that 
\begin{equation}
\zeta' (x+yi) \zeta (x-yi) = - \zeta (x+yi) \zeta' (x-yi) \implies \frac{ \zeta' (x+yi)}{\zeta (x+yi)} = - \frac{ \zeta' (x-yi)}{\zeta (x-yi)}. 
\end{equation}
Integrating the above expression we find 
\begin{equation}
\ln{\left(\zeta (x+yi) \right)} = \ln{\left( \left( \zeta (x-yi) \right)^{-1} \right)} \implies \zeta (x+yi) = \frac{1}{ \zeta (x-yi)} \implies  \zeta (x+yi) \zeta (x-yi) = 1. 
\label{cond1}
\end{equation}

We can now do the same but using the first-order partial derivative with respect to $y$:
\begin{eqnarray}
\partial_y |\zeta (x+yi)|^2 &=& \partial_y \left(\zeta (x+yi) \zeta (x-yi)\right) = \left( \partial_y \zeta (x+yi) \right) \zeta (x-yi) +  \zeta (x+yi) \left( \partial_y \zeta (x-yi) \right), \nonumber \\ 
&=& i \cdot \zeta' (x+yi) \zeta (x-yi) -i \cdot  \zeta (x+yi)  \zeta' (x-yi) = 0. 
\end{eqnarray}
This equation implies that 
\begin{equation}
\zeta' (x+yi)  \zeta (x-yi) = \zeta (x+yi)  \zeta' (x-yi)  \implies \frac{ \zeta' (x+yi)}{\zeta (x+yi)} =  \frac{\zeta' (x-yi)}{\zeta (x-yi)}. 
\end{equation}
Integrating the above expression we find 
\begin{equation}
\ln{\left(\zeta (x+yi) \right)} = \ln{\left( \zeta (x-yi) \right)} \implies \zeta (x+yi) = \zeta (x-yi). 
\label{cond2}
\end{equation}

Hence, using Eq.\ref{cond1} in combination with Eq.\ref{cond2} we have: 
\begin{equation}
\left( \zeta(x + yi) \right)^2 = 1. 
\end{equation}
Therefore, the question is now then rephrased as if there exists an $x + y i \in \mathbb{C}$ with $0 < a < x < 1-a < 1$ and such that $\zeta(x + yi) = +1$ or $-1$. 

Notice that if $\zeta(x + yi) = \pm 1$, then ${\rm Re}\left(\zeta(x + yi)\right) = \pm 1$ and ${\rm Im}\left(\zeta(x + yi)\right) = 0$, implying then that $|\zeta(x + yi) \mp 1| = 0$. In what follows we tackle whether this can be true using several methods. 

First, we analyzed the above condition numerically by means of extensive Matlab simulations, implementing $\zeta(x + yi)$ both in terms of its series definition and in terms of the functional equation, seeing no difference between the two implementations. We have observed that the function is symmetric with respect to $y$, and therefore studied it only for $y \ge 0$. In Fig.(\ref{figNew}) we show surface plots of $|\zeta(s) \mp 1|$ (with $s = x + yi$). We see that the function is always strictly larger than zero in the considered domains.  Both figures also give an intuition of what is the structure of the respective functions, with oscillations coming from the sinusoidal part in the functional definition of the zeta function. 
\begin{figure}[H]
\begin{centering}
\includegraphics[scale=0.5]{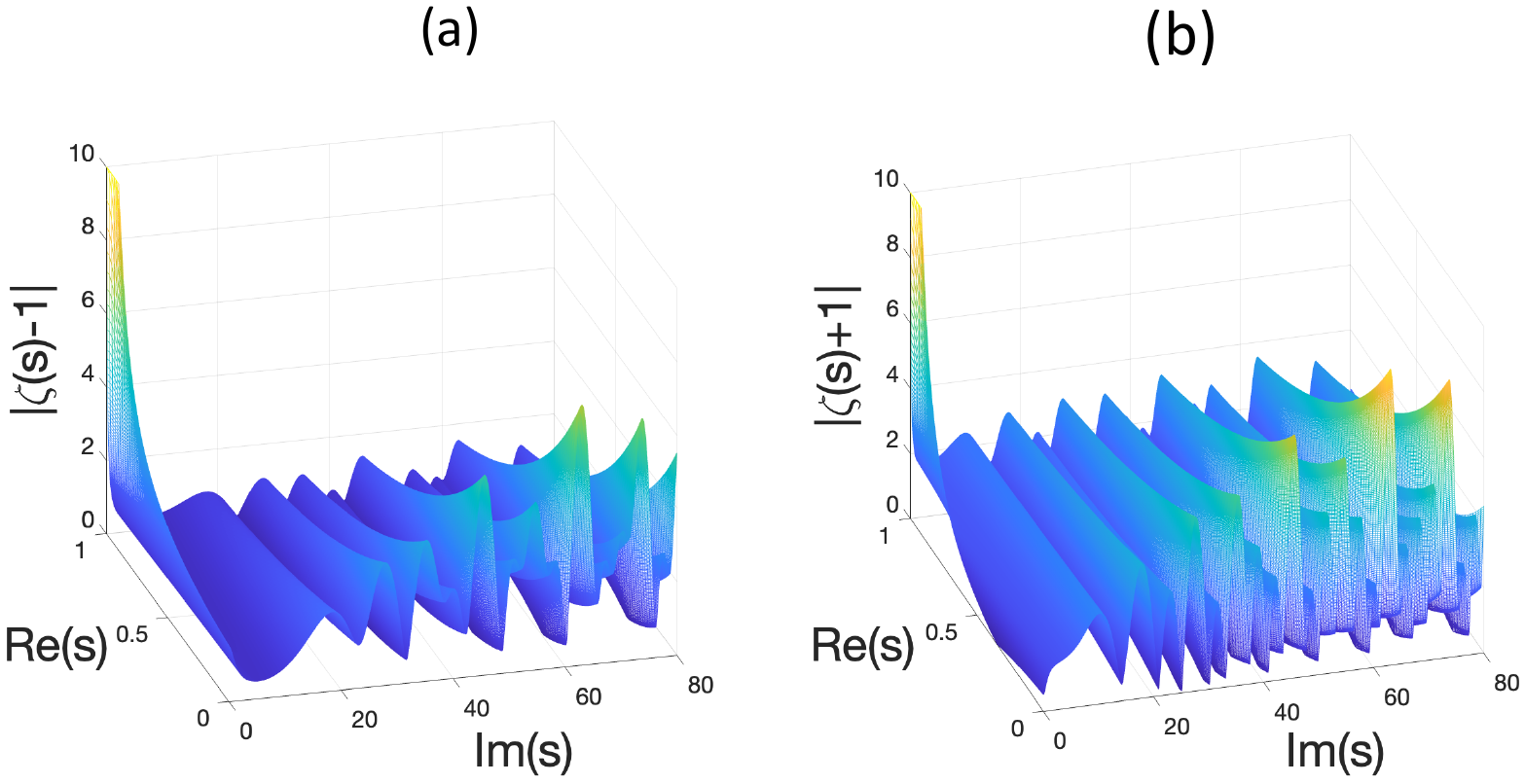}
\par\end{centering}
\caption{Modulus of (a) $\zeta(s) - 1$ and (b) $\zeta(s) + 1$, for $0<{\rm Re}(s)<1$ and $0<{\rm Im}(s)<80$. The function is symmetric with respect to ${\rm Im}(s) \leftrightarrow -{\rm Im}(s)$.\label{figNew}}
\end{figure}

In addition, we have minimized both functions to the best of our numerical capabilities using a gradient-free Nelder-Mead algorithm with Matlab. Up to machine precision, it was impossible to find exact zeros for $0 < x < 1$ (we found the minimum with magnitude around $O(10^{-5})$). The Matlab code is provided in Appendix \ref{apA}. 

Our analysis has also shown that it is possible to find exact zeros for $x \notin (0,1)$. For instance, using Matlab minimization with Nelder-Mead methods we could find exact zeros for $x \approx 1400$. In addition, using WolframAlpha to solve the system of equations ${\rm Re} \left( \zeta(x + yi) \right) = 1$ and ${\rm Im} \left( \zeta(x + yi) \right) = 0$, we found exact zeros at values such as $(x,y) = {(53.7781, 9.0647), (54.8643,-4.5323), (55.2193, -9.0647),...}$, clearly beyond the considered domain, and with similar conclusions for the case ${\rm Re} \left( \zeta(x + yi) \right) = -1$, finding exact zeros at values such as $(x,y) = {(-5.0154, 5.5222), (-5.0154, -5.5222), (-0.3479,14.2831),...}$. Details are also provided in Appendix \ref{apA}. 

Therefore, $|\zeta (x+yi)|^2$ does not have a maximum or a minimum for some $x + yi \in \mathbb{C}$ with $0 < a < x < 1-a < 1$, and the same happens for $|\zeta (x+yi)|$, hence it cannot be possible that $\zeta (a+yi) = \zeta (1-a+yi) = 0$. So, we arrive to a contradiction. 

This implies that if $\zeta(s)$ has any zero with $0<{\rm Re}(s)<1$, then the real part cannot be different to $\displaystyle 1/2$, but instead equal to $\displaystyle 1/2$, as we wanted to prove. This observation is correct within the validity and capacity of our numerical calculations. $\square$ 

}

Next, we consider the second direction of the implication. 

\noindent
\underline{\textit{(ii) If ${\rm Re}(s)=\displaystyle 1/2$ then $\exists y \in \mathbb R$ such that $\zeta(\displaystyle 1/2+yi) = 0$}}

To prove this statement, we should prove that $\exists y\in\mathbb R$ such that $\zeta (\displaystyle 1/2+yi)=0$, that is, such that ${\rm Re}(\zeta (\displaystyle 1/2+yi))=0$ and ${\rm Im}(\zeta (\displaystyle 1/2+yi))=0$.

We use again reduction to absurdity. That is, we will assume that $\nexists y \in\mathbb R$ such that ${\rm Re}(\zeta(\displaystyle 1/2+yi))=0$ and ${\rm Im}(\zeta(\displaystyle 1/2+yi))=0$. And this means that, applying logic rules for the negation of quantifiers and the Morgan law to negate a conjunction, we must prove that $\forall y\in \mathbb R$ it is satisfied that ${\rm Re}(\zeta(\displaystyle 1/2+yi)) \neq 0$ or ${\rm Im}(\zeta(\displaystyle 1/2+yi)) \neq 0$.

Hence, let us assume that $\forall y\in \mathbb R$ it is satisfied that $\zeta(\displaystyle 1/2+yi)=a+bi \neq 0$. Below we analyze if it can be satisfied that $b\neq 0$ or $a\neq 0$.

\begin{itemize}
	\item $b\neq 0$: In such case this implies that $\forall y\in \mathbb R$ the imaginary part has always the same sign. However, using Proposition 2, $\zeta(\displaystyle 1/2 +(-y)i)=\overline{\zeta(\displaystyle 1/2 +yi)}=\overline{a+bi}=a-bi$. That is, $\forall y\in \mathbb R$, we have that ${\rm Im}(\zeta(\displaystyle 1/2  +yi))$ changes its sign between $-y$ and $y$. Hence, using Bolzano's theorem, it will be zero between $-y$ and $y$. Therefore, we arrive to a contradiction, and this implies that $\exists y \in \mathbb R$ such that ${\rm Im}(\zeta(\displaystyle 1/2 +yi))=b=0$.
    \item $a\neq 0$: In such case this implies that $\forall y\in \mathbb R$ the real part has always the same sign. However, we know that $\zeta(\displaystyle 1/2)= -1.4603545088...$, and for instance, $\zeta(\displaystyle 1/2+5i)= 0.7018123711... + 0.2310380083... i$. That is, for $y=0$ the real part is negative, whereas for $y=5$ it is positive. Hence, using Bolzano's theorem, the real part will be zero between $y=0$ and $y=5$. Therefore, we arrive again to a contradiction, and this implies that $\exists y \in \mathbb R$ such that ${\rm Re}(\zeta(\displaystyle 1/2 +yi))=a=0$. The plot in Fig.\ref{fig7} shows this case.

\begin{figure}[H]
\begin{centering}
\includegraphics[scale=0.4]{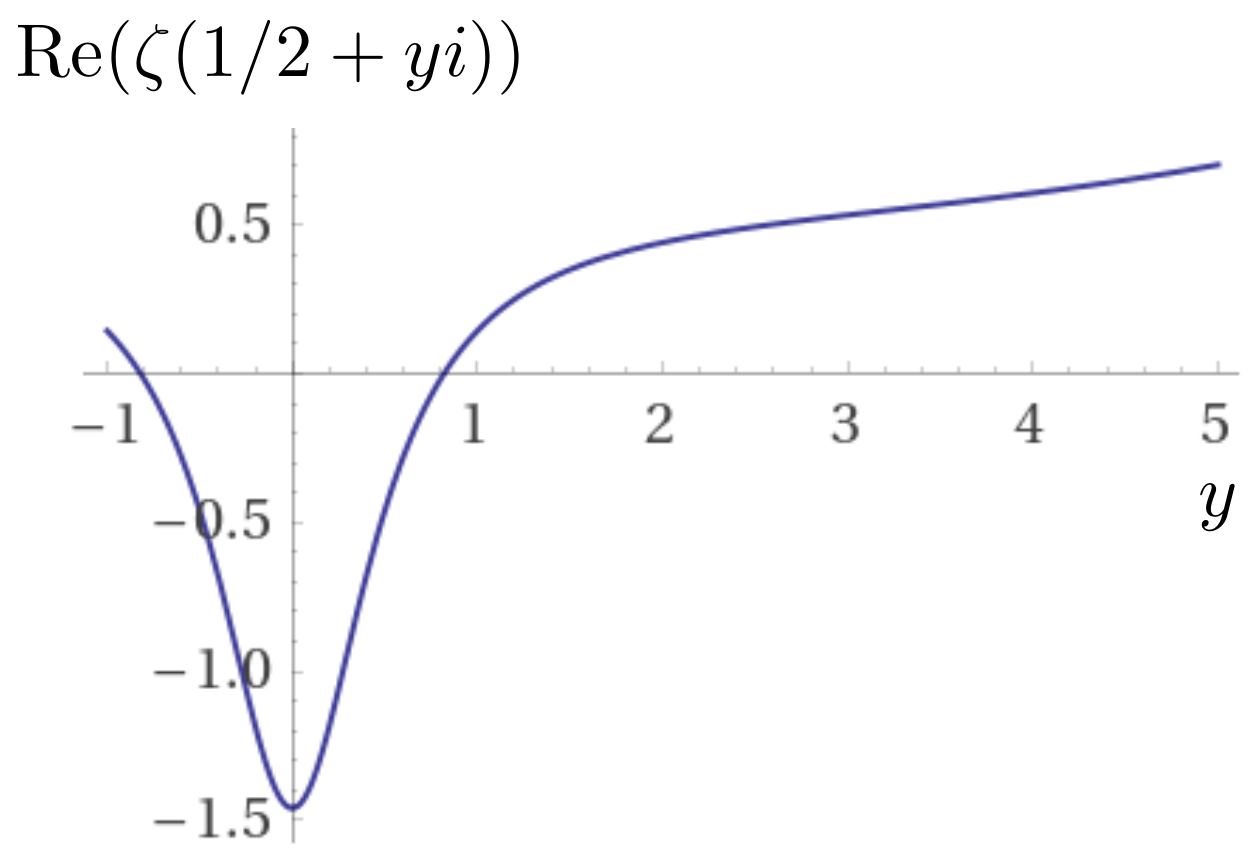}
\par\end{centering}
\caption{Real part of $\zeta(\displaystyle 1/2+yi)$ for $-1<y<5$. \label{fig7}}
\end{figure}
\end{itemize}

Therefore, in both cases we arrive to a formal contradiction, and this implies that it cannot be possible that  $\nexists y \in \mathbb R$ such that ${\rm Re}(\zeta(\displaystyle 1/2+yi))=0$ and ${\rm Im}(\zeta(\displaystyle 1/2+yi))=0$, or what is the same, that $\forall y \in \mathbb R$ it is satisfied that ${\rm Re}(\zeta(\displaystyle 1/2+yi)) \neq 0$ or ${\rm Im}(\zeta(\displaystyle 1/2+yi)) \neq 0$. 

So, this implies that $\exists y \in \mathbb R$ such that $\zeta (\displaystyle  1/2+yi)=0$, as we wanted to prove. $\square$

\section{Constraints on the zeros}
\label{seczeros}

Applying small adjustments to the above analysis, our next question is: how are these $y \in \mathbb R$ such that $\zeta (\displaystyle  1/2+yi)=0$?. Let us analyze this below.

Let $y \in \mathbb R$ be a real number such that $\zeta (\displaystyle 1/2+yi)=a+bi \neq 0$, that is, such that the real and imaginary parts of $\zeta (\displaystyle 1/2+yi)$ are not simultaneously $0$. Then, using the functional equation Eq.(\ref{eq1}) we have that 
\begin{equation} \label{eq20}
\zeta \left(\displaystyle \frac {1}{2}+yi \right)=a+bi=\sqrt{\displaystyle \frac {2}{\pi}}(2\pi)^{yi}\sin{\left(\displaystyle \frac{\pi}{4} + \displaystyle \frac {y\pi}{2}i\right)}\Gamma\left(\displaystyle \frac {1}{2} - yi\right)\zeta \left(\displaystyle \frac {1}{2}-yi\right). 
\end{equation}
And using Proposition 2 we know that $\zeta(\displaystyle 1/2-yi)=a-bi$. So, this implies that 
\begin{equation} \label{eq21}
a+bi=\sqrt{\displaystyle \frac {2}{\pi}}(2\pi)^{yi}\sin{\left(\displaystyle \frac{\pi}{4} + \displaystyle \frac {y\pi}{2}i \right)}\Gamma\left(\displaystyle \frac {1}{2} - yi \right)(a-bi). 
\end{equation}

Then, since $a-bi \neq 0$, we have that 
\begin{eqnarray}
\displaystyle \frac{a+bi}{a-bi}&=&\sqrt{\displaystyle \frac {2}{\pi}}(2\pi)^{yi}\sin{\left(\displaystyle \frac{\pi}{4} + \displaystyle \frac {y\pi}{2}i \right)}\Gamma\left(\displaystyle \frac {1}{2} - yi \right), \nonumber \\
\displaystyle \frac {(a+bi)^{2}}{(a-bi)(a+bi)}&=&\sqrt{\displaystyle \frac {2}{\pi}}(2\pi)^{yi}\sin{\left(\displaystyle \frac{\pi}{4} + \displaystyle \frac {y\pi}{2}i \right)}\Gamma\left(\displaystyle \frac {1}{2} - yi \right), \nonumber \\
\displaystyle \frac {a^{2}-b^{2}+2abi}{a^{2}+b^{2}}&=&\sqrt{\displaystyle \frac {2}{\pi}}(2\pi)^{yi}\sin{\left(\displaystyle \frac{\pi}{4} + \displaystyle \frac {y\pi}{2}i \right)}\Gamma\left(\displaystyle \frac {1}{2} - yi \right), \nonumber \\
\displaystyle \frac {a^{2}-b^{2}}{a^{2}+b^{2}} + \displaystyle \frac {2ab}{a^{2}+b^{2}}i&=&\sqrt{\displaystyle \frac {2}{\pi}}(2\pi)^{yi}\sin{\left(\displaystyle \frac{\pi}{4} + \displaystyle \frac {y\pi}{2}i \right)}\Gamma\left(\displaystyle \frac {1}{2} - yi \right). 
\end{eqnarray}
Hence, 
\begin{eqnarray} \label{eq23}
{\rm Re}\left(\sqrt{\displaystyle \frac {2}{\pi}}(2\pi)^{yi}\sin{\left(\displaystyle \frac{\pi}{4} + \displaystyle \frac {y\pi}{2}i\right)}\Gamma\left(\displaystyle \frac {1}{2} - yi\right)\right)&=&\displaystyle \frac {a^{2}-b^{2}}{a^{2}+b^{2}} ,\nonumber \\
{\rm Im}\left(\sqrt{\displaystyle \frac {2}{\pi}}(2\pi)^{yi}\sin{\left(\displaystyle \frac{\pi}{4} + \displaystyle \frac {y\pi}{2}i\right)}\Gamma\left(\displaystyle \frac {1}{2} - yi\right)\right)&=&\displaystyle \frac {2ab}{a^{2}+b^{2}}. 
\end{eqnarray}

\begin{figure}[H]
\begin{centering}
\includegraphics[scale=0.5]{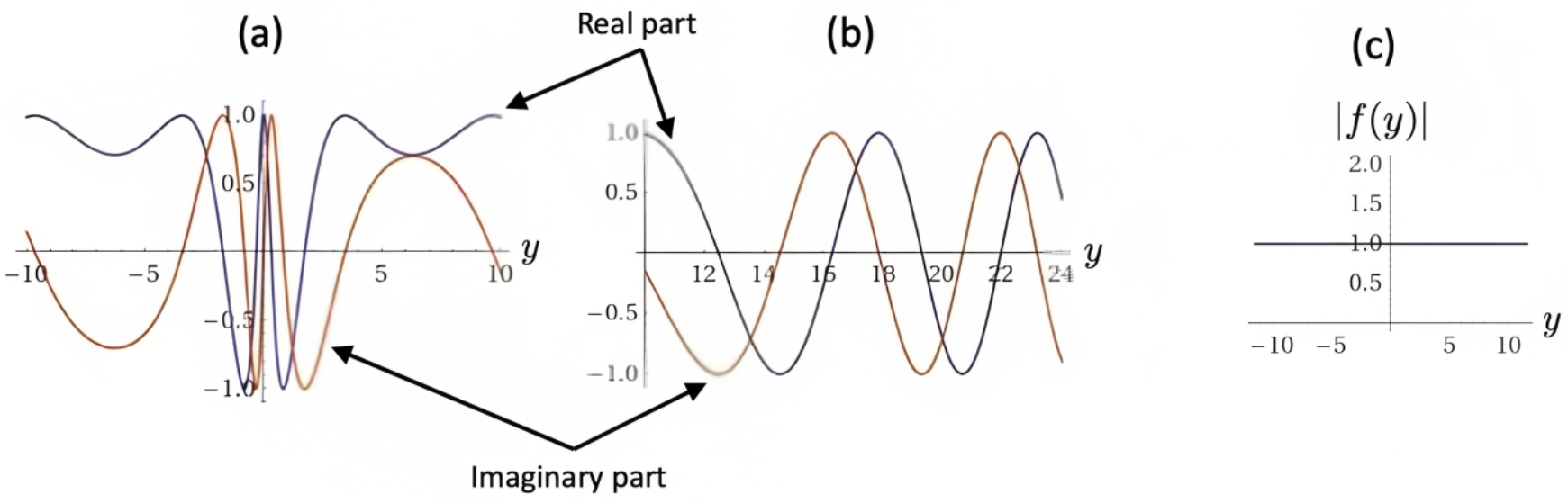}
\par\end{centering}
\caption{Function $f(y)$ for (a) $-10<y<10$ and (b) $10<y<24$, separated for ease of identification of different behaviour regimes. (c) Modulus $|f(y)|$ for $-10 < y < 10$. The same behaviour for the modulus is observed no matter the value of $y$. \label{fig8}}
\end{figure}

Let us now call 
\begin{equation}
f(y) = \sqrt{\displaystyle \frac {2}{\pi}}(2\pi)^{yi}\sin{\left(\displaystyle \frac{\pi}{4} + \displaystyle \frac {y\pi}{2}i\right)}\Gamma\left(\displaystyle \frac {1}{2} - yi\right). 
\label{neweq}
\end{equation}
The behavior of this function is shown in Figs.\ref{fig8}a,b. We see also in Fig.\ref{fig8}c that the function also satisfies $|f(y)|=1$ for all $y$. 

\noindent
Then, if we call $c(y)={\rm Re}(f(y))$ and $d(y)={\rm Im}(f(y))$, since $|f(y)|=1$, we have that 
\begin{equation}
(c(y))^{2} + (d(y))^{2}=1. 
\end{equation}
Therefore, if $c(y) =\pm 1$ then $d(y)=0$, and if $d(y) =\pm 1$ then $c(y)=0$. If we now analyze these cases using Eq.\ref{eq23}, we have the following: 

\noindent
\textbf{a) Case $c(y) =- 1$ and $d(y)=0$:}
\begin{eqnarray}
-1&=&\displaystyle \frac {a^{2}-b^{2}}{a^{2}+b^{2}}, \nonumber \\
0&=&\displaystyle \frac {2ab}{a^{2}+b^{2}}. 
\end{eqnarray}
That is:
\begin{eqnarray}
-a^2-b^2 &=& a^2-b^2 , \nonumber \\
0&=&2ab.
\end{eqnarray}

From the first equation we obtain $a=0$, and then the second equation is always satisfied for all $b$, included $b=0$. Therefore, if $c(y) =- 1$ and $d(y)=0$, then $a+bi$ can be equal to $0$.

\noindent
\textbf{b) Case $c(y) = 1$ and $d(y)=0$:}
\begin{eqnarray}
1&=&\displaystyle \frac {a^{2}-b^{2}}{a^{2}+b^{2}}, \nonumber \\
0&=&\displaystyle \frac {2ab}{a^{2}+b^{2}}. 
\end{eqnarray}
That is:
\begin{eqnarray}
&a^2+b^2 = a^2-b^2 ,\nonumber \\
&0=2ab. 
\end{eqnarray}

From the first equation we obtain $b=0$, and then the second equation is always satisfied for all $a$, included $a=0$. Therefore, if $c(y) = 1$ and $d(y)=0$, then $a+bi$ can be equal to $0$.

\noindent
\textbf{c) Case $c(y) = 0$ and $d(y)=1$:}
\begin{eqnarray}
0&=&\displaystyle \frac {a^{2}-b^{2}}{a^{2}+b^{2}}, \nonumber \\
1&=&\displaystyle \frac {2ab}{a^{2}+b^{2}}. 
\end{eqnarray}
That is:
\begin{eqnarray}
0&=& a^2-b^2 ,\nonumber \\
a^2+b^2&=&2ab. 
\end{eqnarray}

From the second equation we obtain $(a-b)^{2}=0$, ie $a=b$, and then the first equation is always satisfied, included for $a=b=0$. Therefore, if $c(y) = 0$ and $d(y)=1$, then $a+bi$ can be equal to $0$.

\noindent
\textbf{d) Case $c(y) = 0$ and $d(y)=-1$:}
\begin{eqnarray}
0&=&\displaystyle \frac {a^{2}-b^{2}}{a^{2}+b^{2}} ,\nonumber \\
-1&=&\displaystyle \frac {2ab}{a^{2}+b^{2}}.  
\end{eqnarray}
That is:
\begin{eqnarray}
0&=&a^2-b^2, \nonumber \\
-a^2-b^2&=&2ab. 
\end{eqnarray}

From the second equation we obtain $(a+b)^{2}=0$, ie $a=-b$, and then the first equation is always satisfied, included for $a=b=0$. Therefore, if $c(y) = 0$ and $d(y)=-1$, then $a+bi$ can be equal to $0$.

Hence, as a consequence of the results obtained above, we can not assume that $\zeta (\displaystyle 1/2+yi)=a+bi \neq 0$ for all $y \in \mathbb R$. Therefore, $\exists y \in \mathbb R$, satisfying any of the four cases analyzed above, such that $\zeta (\displaystyle 1/2+yi)=a+bi=0$, as we wanted to prove. $\square$

\section{Conclusions}
\label{sec4}

In this paper we have analyzed the Riemann hypothesis. \hlt{Our objective has been, from the beginning, to validate whether this hypothesis is true or not. 

We did such validation using the functional equation of the zeta function, and stating that, if and only if the real part of a complex number was equal to $1/2$, only then could we find a real number $y$ such that $\zeta(1/2+iy)$ was equal to zero, so that the non-trivial zeros of the function all lie on the strip line $1/2+iy$.

We began our study by stating and proving two propositions. The first proposition was \emph{the key step} in order to start building the whole analysis that we presented in this paper. 

For this analysis, we have used} the functional equation of $\zeta (s)$ for $s \in \mathbb C $ with $0<{\rm Re} (s)<1$, and using the absurdity reduction method, we assumed that the zeros of $\zeta (s)$ can have their real part equal to $a \neq \displaystyle 1/2$ and $0<a<1$, reaching a contradiction. This is obtained after an analytical study based on a complex function, and its modulus as a real function of two real variables, combined with an intensive numerical analysis up to the best of our computational resources in order to check the strict positivity of the given modulus function. Additionally, using also the absurdity reduction method, in combination with some logic rules to negate a quantifier and a Morgan's law, we showed analytically that if ${\rm Re}(s)=\displaystyle 1/2$, then $\exists y \in \mathbb R$ such that $\zeta(\displaystyle 1/2 +yi)=0$. Moreover, we provided analytical conditions that should be satisfied by the $y$ values candidates, so that $\zeta(\displaystyle 1/2 +yi)=0$.

\hlt{We need to restate that the numerical analysis methods implemented with Matlab and Wolfram Alpha work under certain typical assumptions, for example not exceeding a certain accepted minimum and small enough error, also limitations coming from potential sources of error such as rounding in floating point operations and machine precission.}

This is another work toward assessing the validity of Riemman's conjecture. \rom{As such, we have proven here that the conjecture is true, up to the best of our numerical analysis for the strict positivity of $|\zeta{(x + yi)} \mp 1|$ for $0 < x < 1$. Our analysis shows that this function is never exactly zero in this interval,} to the best of our computational power, and we see no tendency for a change. Given that the rest of the derivations here are fully analytical, we conclude that Riemann's conjecture is true. 

\bigskip 

\noindent {{\bf Acknowledgements:}} \rom{We acknowledge Prof. Richard Taylor, from Stanford University, for an insightful comment on a previous version of this manuscript.} Fruitful discussions over the years with many mathematicians and physicists' colleagues about the validity and applications of the conjecture are also acknowledged. Finally, we also acknowledge computational support from Mathworks and Wolfram.    

%\section{References}

\appendix 

\section{Appendix: numerical approaches}
\label{apA}

\subsection{Matlab code and Nelder-Mead method} 

What follows is the Matlab code to minimize $|\zeta(s) \mp 1|$: 

\begin{lstlisting}

% Computes the minimum of (abs(zeta(s))-+1) with gradient-free Nelder-Mead
% simplex method. It takes many possible initial conditions for the
% optimization, and evaluates the best solution over all of them. 

clear all 

x = 1; 
for a = 0:0.1:1 % real part 
    y = 1; 
    for b = 0:0.1:80 % imaginary part
        fun = @(q)abs(zeta(q(1) + 1i*q(2))-1); % replace "-1" by "+1" for the second case. For the functional form (no numerical difference), use instead fun = @(q)(abs((2^(q(1) + 1i*q(2)))*(pi^((q(1) + 1i*q(2))-1))*sin(pi*(q(1) + 1i*q(2))/2)*gamma(1-(q(1) + 1i*q(2)))*zeta(1-(q(1) + 1i*q(2)))+1)). 
        q0 = [a,b]; % initial condition of optimization
        [q,fval] = fminsearch(fun,q0); % Nelder-Mead simplex method 
        minfun(x,y) = fval;
        minxx(x) = q(1); 
        minyy(y) = q(2); 
        dat = {minxx minyy minfun};
        save('Minus.mat', 'dat') % save data
        [a b q(1) q(2) fval] % data on screen
        y = y+1; 
    end
    x = x+1; 
end

figure; mesh(minfun) % Plots distribution of found minima. 

minimum = min(min(minfun));
[rr,ss]=find(minfun==minimum); 
[minxx(rr), minyy(ss), minimum] % Lowest-found minima. 

\end{lstlisting}

Matlab allows the implementation of different optimization functions, including gradient-based and gradient-free methods. The code above uses a gradient-free simplex method called Nelder-Mead. This method is particularly suitable for functions such as the one being analyzed here, which has a very chaotic gradient. the Nelder-Mead method (also known as downhill simplex method, amoeba method, or polytope method), is a numerical method used to find extremes of a real function in a multidimensional space. The method uses simplices and approximates a local optimum of a problem. It does so by maintaining a set of test points arranged as a simplex, and then extrapolates the behavior of the objective function measured at each test point in order to find a new test point and to replace one of the old test points with the new one. As such, it is a direct search method  and is often applied to nonlinear optimization problems for which derivatives may not be known.

\subsection{Wolfram Alpha calculations}

The numerical calculations done with Wolfram Alpha were based on solving a system of two equations, where we provide the system with the following constraints:
\begin{equation} 
{\rm Re}(\zeta(x+yi))=\mp 1, ~~~ {\rm Im}(\zeta(x+yi))=0. 
\end{equation}
Wolfram provides flexible tools for numerical root finding using a variety of algorithms, such as Newton's method (which is gradient-based) as well as the bisection method (which is gradient-free). In Newton's method, the idea is to start with an initial guess, then to approximate the function by its tangent line, and finally to compute the intercept of this tangent line with the axis. This point will typically be a better approximation to the original function's root than the first guess, and the method can be iterated. In addition, the bisection method applies to any continuous function for which one knows two values with opposite signs. The method consists of repeatedly bisecting the interval defined by these values, and then selecting the subinterval in which the function changes sign, therefore pinpointing the root. Both methods, Newton and bisection, use as stopping criterion that the difference between the computed points in two consecutive iterations is sufficiently small, i.e., $|x_{k+1} - x_{k}| < \epsilon$, with $x_k$ the point at step $k$ (and similarly for $x_{k+1}$), and $\epsilon$ some error tolerance. In some cases the stopping criteria can be applied also to the relative difference, i.e., $|(x_{k+1} - x_{k})/x_k| < \epsilon$.
{}

\end{document}